\documentclass[conference,compsoc]{IEEEtran}
\ifCLASSOPTIONcompsoc
  \usepackage[nocompress]{cite}
\else
  \usepackage{cite}
\fi
\ifCLASSINFOpdf
\else
\fi
\usepackage{amsmath}
\usepackage{algorithmic}

\usepackage{stfloats}
\usepackage{url}

\usepackage{color}
\usepackage{mathtools}
\usepackage{amsfonts}
\usepackage{amsthm}

\newcommand{\GCD}{\texttt{GCD}}
\newcommand{\R}{\ensuremath{\mathbb{R}}}
\newcommand{\C}{\ensuremath{\mathbb{C}}}

\usepackage{siunitx}
\usepackage{booktabs}
\usepackage{amsmath} 
\usepackage{color}
\usepackage{mathtools}
\usepackage{amsfonts}
\usepackage{amsthm}
\usepackage{float}

\usepackage{csquotes}

\usepackage{listings}
\theoremstyle{definition}
\newtheorem*{definition}{Definition}

\renewcommand{\emph}[1]{\textsl{#1}}

\newtheorem{theorem}{Theorem}

\usepackage{booktabs} 
\usepackage{tikz}
\usepackage{mathdots}
\usepackage{esvect}
\usepackage{bbm}
\usepackage{amsmath}

\usepackage{subfig}

\usepackage{maplestd2e}

\usepackage{xcolor}
\definecolor{darkgreen}{rgb}{0,.35,0}
\definecolor{darkblue}{rgb}{0,0,.5}
\definecolor{darkred}{rgb}{.6,0,0}

\numberwithin{equation}{section}

\usepackage{amsmath}
\usepackage{amssymb}
\usepackage{amsfonts}
\usepackage{graphicx}
\usepackage{hyperref}
\hypersetup{
    linkcolor=blue,
    filecolor=blue,      
    urlcolor=blue,
}

\usepackage{algorithmic}
\usepackage{algorithm}
\usepackage{tikz-cd}
\usepackage{chngcntr}

\usepackage{pgf, tikz}
\usetikzlibrary{arrows, automata}
\makeatletter
\def\Ddots{\mathinner{\mkern1mu\raise\p@
\vbox{\kern7\p@\hbox{.}}\mkern2mu
\raise4\p@\hbox{.}\mkern2mu\raise7\p@\hbox{.}\mkern1mu}}
\makeatother

\newcommand{\agcd}[2]{\ensuremath{\mathrm{agcd}_{\rho}^{\sigma}(#1,#2)}}

\usepackage{xcolor}

\usepackage[T1]{fontenc}
\usepackage{enumitem}

\newcommand\dm[1]{\mathbf{d}_{#1}}

\begin{document}
%
\title{Approximate GCD in Lagrange bases}

\author{\IEEEauthorblockN{Leili Rafiee Sevyeri and Robert M.~Corless }
\IEEEauthorblockA{Ontario Research Centre for Computer Algebra,
Western University, Canada and\\
The David R.~Cheriton School of Computer Science\\
University of Waterloo, 
Canada}
}


%


\maketitle

\begin{abstract}
For a pair of polynomials with real or complex coefficients, given in any particular basis, the problem of finding
their \GCD\  is known to be ill-posed. An answer is still desired for many applications, however. Hence, looking for a 
\GCD\  of so-called \textsl{approximate polynomials} where this term explicitly denotes small uncertainties in the coefficients has received significant attention in the field of hybrid symbolic-numeric computation. 
In this paper we give an algorithm, based on one of Victor Ya.~Pan, to find an approximate \GCD\  for a 
pair of approximate polynomials given in a Lagrange basis. More precisely, we suppose that these polynomials 
are given by their approximate values at distinct known points. We first find each of their roots by using a Lagrange basis companion matrix for each polynomial, cluster the roots of each polynomial to identify multiple roots, and then ``marry'' the two polynomials to find their \texttt{GCD}. At no point do we change to the monomial basis, thus preserving the good conditioning properties of the original Lagrange basis.  We discuss advantages and drawbacks of this method.  The computational cost is dominated by the rootfinding step; unless special-purpose eigenvalue algorithms are used, the cost is cubic in the degrees of the polynomials.  In principle, this cost could be reduced but we do not do so here.
\end{abstract}


%
\IEEEpeerreviewmaketitle

\section{Introduction}
The Euclidean algorithm or its variant, the extended Euclidean algorithm, is one of the most
well-known and useful algorithms in symbolic computation. One important role of this rational algorithm for finding
greatest common divisors of polynomials in computer algebra is to identify, in the case where all polynomials are exact, polynomials that have multiple roots, and indeed to compute a square-free factoring of a polynomial by starting with the computation of the \GCD\ of the polynomial and its derivative. Although computing the
\GCD\  of two polynomials is straightforward in exact arithmetic and although much effort has been expended to find more efficient algorithms than the Euclidean algorithm such as using subresultants, over $\R$ or $\C$
the \GCD\  problem is ill-posed. 

In order to deal with ill-posedness of the \GCD\  problem, one can look at the approximate \GCD\ 
problem. For a given pair $P$ and $Q$ of univariate polynomials over reals, find another
pair $\tilde{P}$ and $\tilde{Q}$ such that they are respectively ``close'' to $P$ and $Q$, but with a nontrivial \GCD.
Their exact \GCD\ is called an \emph{approximate \GCD} of $P$ and $Q$. 
We will make this precise in Section~\ref{sec:definitions}. 

There is a vast literature on the approximate \GCD\  problem, most of which is concerned with polynomials given in power basis. 
There are many interesting results in approximate \GCD\ including but not limited to~\cite{LabahnBeckermann1998, LabahnBeckermannMatos2018, BottingGiesbrechtMay2005, FaroukiGoodman1996, KaltofenYangZhi2007, KarmarkarLkshman1996, KarmarkarLakshman1998, NakatsukasaSeteTrefethen2018, Zeng2004, Zeng2011}
for approximate \GCD\  problem in the power basis. 

Although the power basis is the most well understood and common basis to represent 
polynomials, it is not the only one. For various purposes, sometimes we would rather 
to have a polynomial in another basis such as a Bernstein basis when one wants to find intersection of B\'ezier curves; this in turn has applications
to Computer Aided Geometric Design. A theoretical solution would be to convert the given pair of 
polynomials in a Bernstein basis into the power basis. In practice this is not recommended
because the conversion is numerically unstable, exponentially so in the degree~\cite{Beckermann2000}.

In \cite{CorlessLeili2019Maple} an algorithm for finding an approximate \GCD\  in a Bernstein basis, without changing basis, is given
which follows the ideas of~\cite{Pan2001}.
In this paper we also use these ideas.

In many applications, one needs to do computations with polynomials given by their values at certain points; that is, in a Lagrange basis. As with Bernstein bases, conversion is usually unstable. Further, Lagrange bases themselves are often \emph{well}-conditioned, and working with them directly preserves their conditioning~\cite{CorlessWatt2004}.  See also~\cite{carnicer2017optimal}. 

Direct computations 
in Lagrange bases are investigated in \cite{Amiraslani2004,AmiraslaniCorlessGonzalez-VegaShakoori2004,Corless2004,Shakoori2004,CorlessRezvani2011,AmiraslaniCorlessLancaster2009,amiraslaniCorlessGunasingam2018}. Applications of these ideas and similar are investigated in, for instance,~\cite{aruliah2007geometric,butcher2011polynomial,corless2008barycentric,DIAZTOCA2006222}. In this 
paper we present an algorithm which follows the same ideas as in~\cite{Pan2001} and 
\cite{CorlessLeili2019Maple} and solves the approximate \GCD\ problem in Lagrange bases. 

Our main working assumption is that the given data is \emph{near} to data specifying a pair of polynomials with common roots.  We allow the problem solver or user to specify a tolerance to say roughly just \emph{how} near, but the method of this paper cannot work if the given data is $O(1)$ away from polynomials with multiple roots.  To see this, one may zero out both polynomials with an $O(1)$ change in the polynomial values, or otherwise trivially make them identical.

Our algorithm works with the sets of \textsc{roots} of a given pair of polynomials. Assume $P$ and $Q$
are given with $\deg(P) = n$ and $\deg(Q) = m$, $\sigma >0$ and $\rho_t$ is a metric (or pseudometric) on the space of 
polynomials of degree $t$. 

We first compute all roots of $P$ and $Q$ using the
companion pencil method introduced in~\cite{Corless2004}. The second step is to cluster the roots of each polynomial,
i.e. find close roots and replace them by one single root with an appropriate multiplicity. 
The next step is to form a bipartite weighted graph with sets of vertices which are the clustered roots 
of $P$ and $Q$. For $r $ and $s$ with $P(r) = Q(s) = 0$, $\lbrace r,s \rbrace$ is an edge with
weight $$W(\lbrace r, s\rbrace) = \min\lbrace \mathrm{multiplicity}_P(r), \mathrm{multiplicity}_Q(s) 
\rbrace.$$ 
Then we use a \emph{Maximum Weighted Matching (MWM)} algorithm to find a maximum matching, say $M$, and
using $M$ we form approximate \GCD. The last step is to construct $\tilde{P}$ and
$\tilde{Q}$ where $\rho_{n}(P,\tilde{P}), \rho_{m}(Q,\tilde{Q})\leq \sigma$. This step is simpler 
in a Lagrange basis than it is for a Bernstein basis or the power basis, provided that
$\rho_t$ be a well behaved metric or pseudometric. 

This paper is organized as follows. Section~\ref{sec:definitions} gives the necessary
formal definitions for approximate \GCD. A method for computing roots of polynomials in Lagrange basis
is given in section~\ref{companion pencil}. Section~\ref{sec:clustering} introduces a divide and conquer algorithm for clustering roots
of a polynomial, together with an alternative set of heuristics that uses symmetry. The maximum weight matching problem is discussed in Section~\ref{sec:MWM}. Section~\ref{sec:apxpoly} contains details of computing 
an approximate \GCD\  and cofactors of polynomials $\tilde{P}$ and $\tilde{Q}$. Finally we provide numerical experiments using
a Maple implementation of our algorithm.
More examples and more details are given in Chapter~3 of~\cite{Rafiee:2020:Hybrid}.

\section{Definitions}\label{sec:definitions}
In this brief section we provide necessary formal definitions which
lead to a careful definition of an approximate \GCD. The most crucial 
component of the notion of an approximate \GCD\  is ``closeness''. Hence, 
we need to formally present a metric on the set of degree $n$ polynomials. Although a metric specifies the notion of closeness
completely, a {\em pseudometric} does the job as well.
\begin{definition}~\cite{Kelley2017}
A metric for a set $X$ is a function $d$ on the cartesian product $X \times X$ to the non-negative reals such that for all points $x,y$, and $z$ of $X$,
\begin{enumerate}
\item $d(x,y)=d(y,x)$,
\item (triangular inequality) $d(x,y) + d(y,z) \geq d(x,z)$,
\item $d(x,y)=0$ if $x=y$, and \label{(c)}
\item if $d(x,y)=0$, then $x=y$. \label{(d)}
\end{enumerate}
\end{definition}
A pseudometric on a set $X$ is a non-negative real valued function $\phi$ defined on $X \times X$ and satisfies all the properties of a metric except for possibly for property~\ref{(d)} which is called  the distinguishability property. 
On the other hand, since a polynomial can be presented with different data in Lagrange
basis, regular well-known metrics are not 
well-defined.


In order to give a well-defined metric (or pseudometric) we
have to consider invariants of polynomials such as their roots. 

In \cite{CorlessLeili2019Maple} we used the so-called \emph{root semi-metric} to solve the approximate \GCD\  problem in Bernstein bases. In 
this work we use the same idea except that we consider it only on 
polynomials with same degree. This leads to a well-defined pseudometric 
which we call it the \emph{root-pseudometric}.

Assume $\rho$ is metric on $\mathbb{C}^n$, $S_n$ is the group of 
permutations of $\lbrace 1, \ldots, n \rbrace$. Moreover, for $f \in 
\mathbb{C}[x]_n$ (set of univariate polynomials of degree $n$) $\lbrace f_1, \ldots , f_n \rbrace$ is the set of roots of $f
$. Let $R_f = \left( f_1, \ldots , f_n \right)$ and $\tau \in S_n$ acts
on $R_f$ by acting on its coordinates, i.e.
$$\tau(R_f) = \left( f_{\tau(1)}, \ldots , f_{\tau(n)} \right).$$
The map
\begin{equation}\label{pseudometric}
\begin{aligned}
 \dm{n}:& \mathbb{R}[x]_n \times \mathbb{R}[x]_n \; \longrightarrow \; \mathbb{R}_{\geq 0} \\
 & \left(f(x),g(x)\right) \;\;\;\; \longmapsto \;\;\;\; \frac{1}{n}\min_{\substack{\tau \in S_n}}
 \Big\lbrace \rho\left( \tau(R_f) , R_g \right) \Big\rbrace
\end{aligned}
\end{equation}
is a pseudometric on $\mathbb{R}[x]_n$. 
Let $$\mathrm{Rep}_{\dm{n}}(f,g) = \big\lbrace \tau \in S_n: \mathbf{d}_n(f,g)= \frac{1}{n}\rho(\tau(R_f),R_g)\big\rbrace.$$   
In fact the identity element of
 $S_n$ gives $\dm{n}(f,f) = 0$ and 
$$\tau \in \mathrm{Rep}_{\dm{n}}(f,g) \Longrightarrow
\tau^{-1} \in \mathrm{Rep}_{\dm{n}}(g,f),$$ where $\tau^{-1}$ is the inverse of $\tau$ in $S_n$. Hence $\dm{n}(g,f)= \dm{n}(f,g)$. The last property of a pseudometric is the triangle inequality. For any fixed $h \in \mathbb{R}[x]_n$ assume $\lambda \in 
\mathrm{Rep}_{\dm{n}}(f,g)$, $\tau \in 
\mathrm{Rep}_{\dm{n}}(f,h)$, $\gamma \in 
\mathrm{Rep}_{\dm{n}}(g,h)$. Now if $\dm{n}(f,h)+\dm{n}(g,h) <\dm{n}(f,g)$ then
$$\rho(\tau(R_f),R_h)+ \rho(\gamma(R_g),R_h) < \rho(\lambda(R_f),R_g),$$
since $\rho$ is a metric on $\mathbb{C}^n$, 
\begin{small}
\begin{align*}
\rho\left(\tau(R_f),\gamma(R_g)\right) \leq & \rho(\tau(R_f),R_h)+ \rho(\gamma(R_g),R_h) \\
 < & \rho(\lambda(R_f),R_g).
\end{align*}
\end{small}
This in particular shows that $\rho(\lambda(R_f),R_g)$ was not minimal, which is a contradiction.

It is worth mentioning that the distinguishability property does not hold for $\dm{n}$. For a non-zero constant $c$, $\dm{n}(cf,f)=0$.

Now that we have a metric in hand, we can formally define an approximate
\GCD\  for a pair of polynomials. A \emph{pseudogcd} set for the pair $P$ and $Q$  is defined as
\begin{multline}
A_{\alpha, r, \sigma} = \big\lbrace g(x) \;\; | \;\;  \exists \tilde{P}, \tilde{Q}\;\; \text{with}\;\; \dm{\deg(P)} (P - \tilde{P})\leq \sigma, \\
\dm{\deg(Q)} (Q- \tilde{Q}) \leq \sigma \;\; \text{and}\;\; g(x) = \gcd(\tilde{P},\tilde{Q}) \big\rbrace \>.
\end{multline}
An approximate \GCD\  for $P,Q$ which is denoted by $\agcd{P}{Q}$, is $G(x) \in A_{\alpha,r, \sigma}$ where 
$\deg(G) =  \max_{\substack{g \in A_{\alpha,r, \sigma}}} \deg(g(x)),$ and $\dm{\deg(P)} (P - \tilde{P}), \dm{\deg(Q)}( Q - \tilde{Q})$ are minimal.

\section{Finding roots of a polynomial in Lagrange basis} \label{companion pencil}
It is well known that solving polynomial equations, or finding the eigenvalues of matrix polynomials, can be done by converting to a generalized eigenvalue problem. In this section we want to find the roots of a given polynomial in Lagrange basis. There are various methods for solving this problem. For example in~\cite{Fortune2001} the author uses a method of Smith~\cite{Smith1970} that has not been generalized to matrix polynomials. However, in this work we prefer to use the generalized eigenvalue problem to find the roots of a certain polynomial.

We recall the method introduced in~\cite{Corless2004} for finding zeros of a given polynomial in Lagrange basis. Consider $p_0, p_1, \ldots, p_n$ which are the values of a polynomial $P$ at $x=x_0$, $x=x_1, \cdots $, and $x_n$. Then \emph{the generalized companion matrix} of the polynomial given by its value is
\begin{equation}
\textbf{C}_0= \left[ \begin{matrix}
0 & -p_0 &-p_1 &\cdots & -p_n \\
\ell_0& x_0 & & &  \\
\ell_1& & \ddots & & \vdots \\
\vdots& & &  &  \\
 \ell_n&  & \cdots &  & x_n\\
\end{matrix}
\right] 
\end{equation}
where $\ell_k = 1 /\prod_{j \neq k} (x_k - x_j)$ are the (scalar) normalization factors of the Lagrange polynomials $L_k (x ) = \ell_k \prod_{j \neq k}(x - x_j)$ and 
$\textbf{C}_1$ is the identity matrix except the $(1,1)$ entry is replaced by zero.
Then we have
\begin{theorem}~\cite{Corless2004}
With the above notation and assumptions
\begin{equation*}
\det ( x \textbf{C}_1 - \textbf{C}_0) = P(x).
\end{equation*}
\end{theorem}

\section{Clustering the roots}\label{sec:clustering}
As we discussed in the first section, our algorithm gets two polynomials as input. In the previous section we explained the method of which 
our algorithm is using to find roots of each polynomial. We note that our method is numerical and the computed roots are approximations. 
So it is natural to wonder if a root with a multiplicity more than~$1$ is approximated as different roots which are close. Although it is not
so clear when this really has happened, we believe working with roots while merging close enough roots is a reasonable approach and likely to lead to good results. This corresponds to what Kahan terms projecting onto the ``pejorative manifold''.  We call the 
process of identifying close roots as the same root with appropriate multiplicity, \emph{root clustering}.

On the other hand, we can also execute this step in the algorithm after constructing an approximate \GCD. However, to be consistent with the
other results in the literature, we put clustering as the second step of our algorithm.

The problem of properly clustering complex roots of perturbed polynomials is hard, and depends strongly on the model used for the perturbation.  In the papers~\cite{imbach2018implementation,becker2016complexity}, for instance, we find the use of an \emph{oracle}, which assumes that information about the \emph{exact} or ``underlying'' polynomial is available at the
request of the user (which might be a program) and in exchange for more computational cost.  That model has been used in polynomial computation since the pioneering work of Sch\"onhage and leads to significant algorithms of practical interest.

However, our model is different: as in~\cite{CorlessGianniTragerBarryWatt1995}
we assume that no oracle is available, and all we have is noisy data about the polynomial (in the Lagrange case, approximate values at well-specified nodes).  In some sense this makes the problem actually impossible: there may be several pairs of polynomials nearby that have common roots. From the incomplete data that we have, however, we must recover the approximate \GCD\ with our best estimate of the multiplicities.  To do this, we resort to heuristics. We present one such algorithm, that uses symmetry and distance, and we also present a divide and conquer algorithm which follows a different idea. 

\subsection{Heuristics for clustering complex roots}
Our first heuristic is the use of distance.  When a polynomial with a multiple root is perturbed, say to $q(z)(z-r)^m + \Delta p(z)$, then the multiple root is, to first order, if $\Delta p(r) = s$ is small, 
\begin{equation}
    z \approx r + \left( -\frac{\Delta p(r)}{q(r)} \right)^{1/m}\>.
\end{equation}
This leads to a cluster of $m$ simple roots at a distance $(s/|q(r)|)^{1/m}$ away, and if $s$ is ``small'' these will be equally-spaced about a near-circle with that radius.  In our current work we do not estimate the size of $q(r)$ but rather use some modest constants or ``fuzz factors'' for this; in future work we will refine this by using the division algorithm of~\cite{Amiraslani2004} to give an estimate of $q(r)$ itself.

The second heuristic is the use of the equal-angle property mentioned above.  Indeed this seems to be quite a reliable indication that we have the correct multiplicity, although we need more experimentation to assert true confidence.

The third heuristic is to look through the list of approximate roots for clusters of higher multiplicity \emph{first}.  This supplies a bias towards higher multiplicity, or equivalently to projection onto the so-called \emph{pejorative manifold}.

Access to these heuristics is possible only because we are looking at approximate roots as a method of finding approximate
\texttt{GCD}. Such heuristics are likely to be less valuable for other approaches to approximate \GCD\ in the Lagrange basis, e.g. via the B\'ezout matrix~\cite{Shakoori2004,diaz2013using,aruliah2007geometric}.

The heuristics are easily fooled for roots of \emph{very} high multiplicity, so our routine has a default maximum order of multiplicity of $3$.  We have used it successfully on artificial examples of multiplicity as high as $5$, however.

\begin{figure}[H]
    \centering
    \includegraphics[width=8cm]{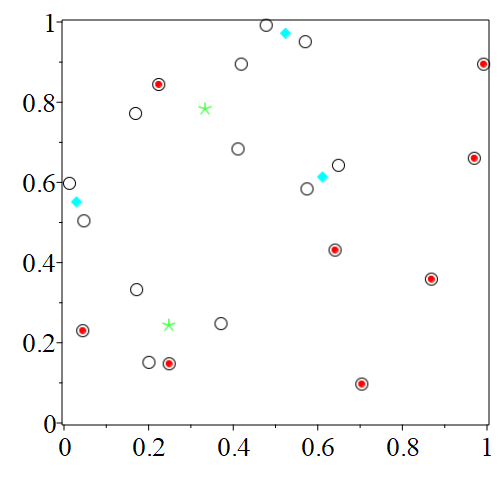}
    \caption{Original data plotted with open black circles.  Clustered triple points plotted with green asterisks.  Clustered double points plotted with cyan diamonds.  Points left untouched by clustering plotted with solid red circles. Note that the triple point clusters each preferred a farther point to include than one that could have been chosen: evidently the heuristic for symmetry is working. }
    \label{fig:clustertwentyrandompoints}
\end{figure}

Let us exhibit our heuristics with a pseudorandom sample of $n=20$ points in the square $[0,1]\times[0,1]$; we used Maple's \texttt{rand()} function twenty times, dividing its integer result by $10^{12}$; half we took for $x$-values and the other half for $y$-values in $z=x+iy$.  We clustered these points with a tolerance of $1/n^2$: one expects a few of $n$ points in the unit square to be $O(1/n^2)$ close.  The results are plotted in Figure~\ref{fig:clustertwentyrandompoints}.

Two sets of triple points were found.  Each of them has the interesting feature that there is one unclustered point closer to the triple's centre than one of the ones chosen.  We believe that the more symmetrical point was chosen, even though it was farther away.

Three sets of double points were found; from the graph, these choices seem very reasonable.
%

\subsection{A Divide and Conquer Clustering Algorithm}

In this subsection we present a deterministic clustering algorithm which
only considers a tolerance directly to cluster a set of given points in 
the complex plain. Our algorithm is a slight modification of the divide and conquer, {\it Closest Pair Algorithm}~\cite[P. 958-960]{CormenLeisersonRivestStein2001}.

Assume $Q$ is an array of  roots of a polynomial. We apply a sorting algorithm to $ Q = (r_1, \ldots , r_n)$ and get the 
array of roots sorted with respect to their real parts. This can be done in $O(n \log n)$. We also add multiplicities of $r_i$'s to them by replacing each of them by $
(r_i,mult(r_i))$. 

\begin{figure}[H]
\includegraphics[scale=.2]{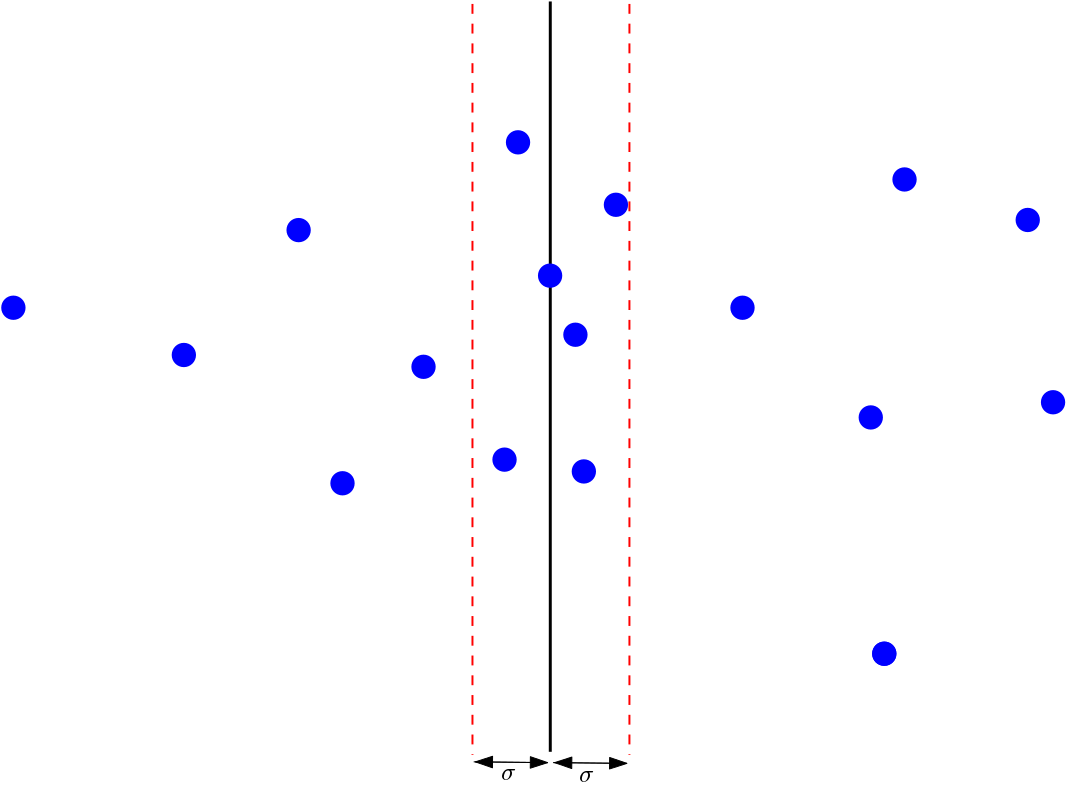}
\centering
\caption{Middle strip for merging the left and right results in clustering 
algorithm}\label{fig:midstrip}
\end{figure}

Our algorithm has $3$ steps:\\
\noindent
{\bf Step $1$:} Split the points into two almost equal parts, $Q_L$ and $Q_R$ by
considering the line $L$, defined as $x = \Re(Q[m][1])$ where $m = \lfloor \frac{\ell+r}{2}\rfloor$, $\ell$ and $r$ are left and right boundaries in recursive calls.
Recursively call the algorithm on $Q_L$ and $Q_R$ then merge the sets (consider imaginary parts for merging).\\
\noindent
{\bf Step $2$:} Once $Q_L$ and $Q_R$ are clustered, form an strip containing 
the points from them which have distance less than $\sigma$ to middle line
(see Figure~\ref{fig:midstrip}). \\
\noindent
\noindent
{\bf Step $3$:} Check the strip for close pairs and merge them. Note that 
for this step $(u,d_u)$ and $(v,d_v)$ will be  merged as 
$$ \left( \dfrac{d_u\cdot u +d_v \cdot v}{d_u+d_v} \right)^{d_u+d_v},$$
if $\vert u-v \vert \leq \sigma$.
%
%
%

\begin{algorithm}
\caption{\texttt{{ClusterRoots}($Q,\ell,r,\sigma$)}}\label{alg:clustercomplex}
\begin{algorithmic}
\REQUIRE A sorted list of roots with their multiplicities, $Q$. \\
\hspace{1cm}$\ell$ is the left boundary (initially $\ell \gets 1$)\\
\hspace{1cm}$r$ is the right boundary (initially $r \gets \texttt{Size}(Q)$). \\
\hspace{1cm}$\sigma$ is the tolerance for clustering.
\ENSURE Clustered roots.
$S \gets Q$\\
{\bf if} $\ell = r$  {\bf then} \\
\hspace{.5cm}$S \gets [Q[l]]$\\
{\bf else} \\
\hspace{.5cm}  $m \gets \lfloor (\ell+r)/2 \rfloor$ \\
\hspace{.5cm}  $S_{L} \gets$ \texttt{ClusterRoots}($Q,\ell, m,\sigma$)\\
\hspace{.5cm}  $S_{R} \gets$ \texttt{ClusterRoots}($Q,m+1, r,\sigma$)\\
\hspace{.5cm}  $S \gets \texttt{Merge}_\Im(Q_L,Q_R$)\\
\hspace{.5cm}  $R \gets \texttt{Select}(S, \Re(S[mid][1]), \sigma)$\\
\hspace{.5cm}  $i \gets \texttt{Search}(S_L, R[first])$\\
\hspace{.5cm}  $j \gets \texttt{Search}(S_R, R[last])$\\
\hspace{.5cm}  $T \gets \texttt{CheckStrip}(R,\sigma)$\\
\hspace{.5cm}  $S \gets \texttt{Merge}(S_L[1..i-1], R, S_R[j+1..last])$\\
{\bf return}$(S)$\\
\end{algorithmic}
\end{algorithm}
An argument which shows that the cost of merging is $O(n)$ is presented in~\cite{CormenLeisersonRivestStein2001}. This means the running time is 
recursively given by $T(n) = 2 T(\frac{n}{2})+O(n)$
and using the Master theorem~\cite{CormenLeisersonRivestStein2001} we can 
resolve it to $T(n) \in O(n\log n)$. Hence the total cost is $O(n \log n)$.

Neither of the above algorithms for clustering roots is perfect. The divide 
and conquer algorithm does not give the highest degrees after clustering.
As an example of its behavior, for the input 
$$[[1, 1], [1.5, 1], [2, 1]],$$
and $\sigma = 0.5$ it returns $$[[1.25, 2], [2, 1]],$$
while we know that $[1.5,3]$ would be a better answer according to 
multiplicity. Moreover \texttt{ClusterRoots} considers the given roots as 
points in the complex plane without considering their properties as 
approximate roots of polynomials.
\\
\\
\section{Maximum Weight Matching problem and approximate \GCD\ }
\label{sec:MWM}
In order to find the \GCD\  of two polynomials, we find their common roots. 
In other words, we look at their roots and find matches between
them. Since it is an exact \GCD\  this matching means equality. In the 
approximate \GCD\  case, matching means close enough with respect to a 
metric.

The above approach in looking at approximate \GCD\  suggests us to use 
matching algorithms from graph theory. In this brief section we
explain the necessary concepts and an algorithm from graph theory.

Assume $G = (V, E, W)$ is a weighted graph. A subset $M \subseteq E$, is 
called a matching if no two edges in $M$ share a common 
vertex. A maximum cardinality matching (MCM) is a matching of 
the maximum size as a set. 

The weight of a matching (in a weighted graph) is defined as
$$W_M = \sum_{e\in M} W(e).$$
A matching $M$ with the maximum weight is called a maximum weight matching 
(MWM).

\begin{figure}
\includegraphics[scale=.5]{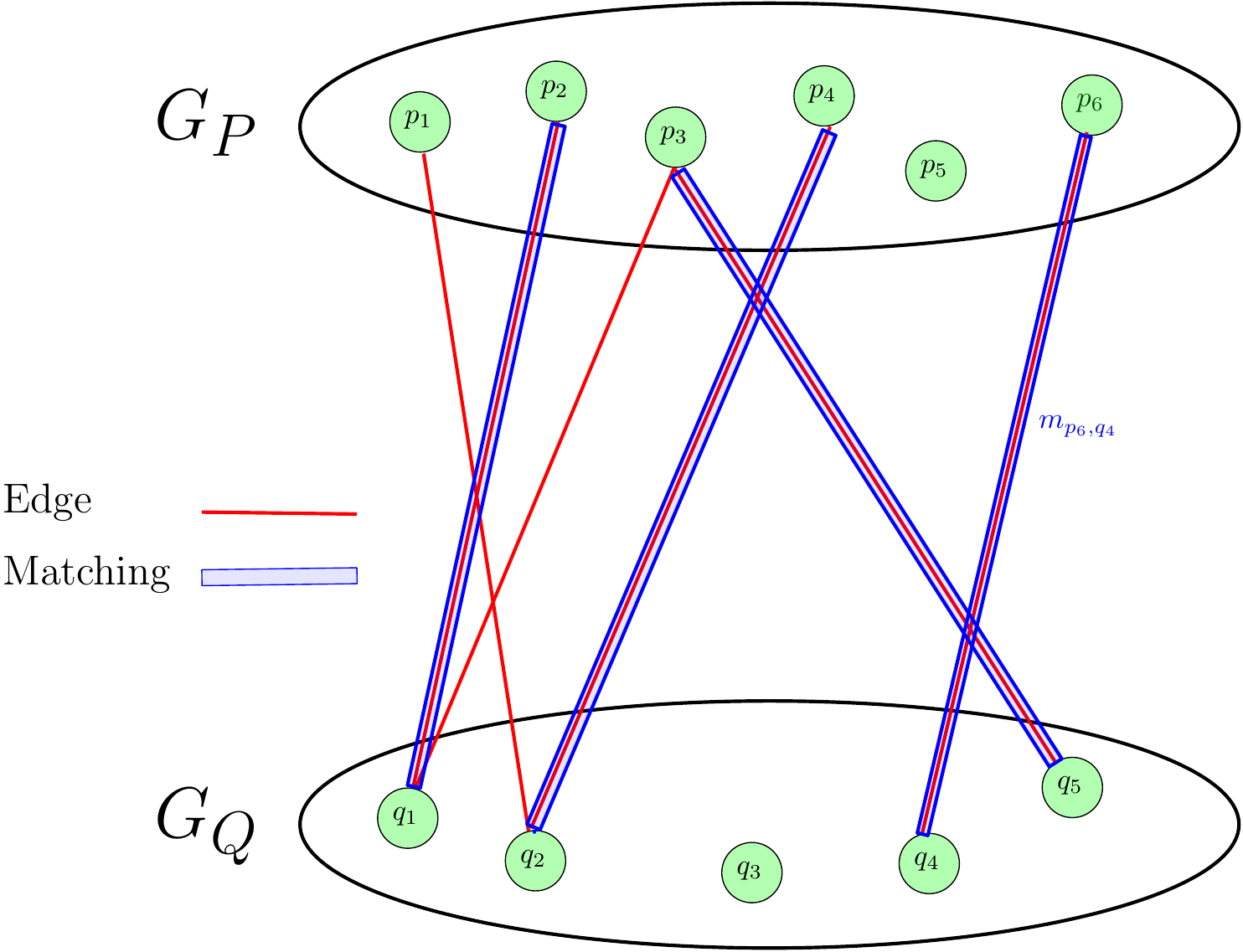}
\centering
\caption{ $m_{p_6,q_4}$ is the minimum of multiplicities of $p_6$ and 
$q_4$.}
\end{figure}
Suppose $n = \vert V \vert$ and $m = \vert E \vert$. It is well known that 
MCM for a bipartite graph can be solved in time $O(m \sqrt{n})$ using the 
deterministic algorithm introduced by Hopcoft and Karp. A more general 
deterministic algorithm for a general graph was given by Micali and 
Vazirani in 1980 which runs in $O(m \sqrt{n})$. There are many more 
algorithms for solving the same problem, we just mention that a randomized
algorithm to find a MCM in a general graph was given by Harvey in 2006 
with complexity $O(n^\omega)$ where $\omega$ is the exponent of 
$n \times n$ matrix multiplication and \cite{LeGall2014} proves that $\omega \leq 2.373$.

In the case that the weights are integers, we define $N$ to be the largest
magnitude of a weight. An MWM for a bipartite graph can be found in $O(N n^\omega)$ using the randomized algorithm introduced by Sankowski in $2006$.

Following~\cite{DuanPettie2014}, for $\delta \in [0,1]$ we call a matching $\delta$-approximate matching if its weight is at least a factor $\delta$ 
of the optimal matching. Similarly a $\delta$-MWM is the problem of finding  a $\delta$-approximate maximum weight matching. Duan and Pettie introduced
an algorithm to solve $\delta$-MWM for a graph with integer weights in
$O(m(1-\delta)^{-1}\log N)$.

A much simpler approximate MWM algorithm is a greedy algorithm. It is well-known that a greedy algorithm solves the problem for a general graph for $\delta = \frac{1}{2}$ (See~\cite{DuanPettie2014}). In our special case, the graph is bipartite and a very simple one. Our Maple implementation of this greedy algorithm gives very good results in practice. We present the pseudocode of the mentioned greedy algorithm here.
 
\begin{algorithm}[H]
\caption{\texttt{GreedyMWM}($G$)}
	\label{alg:greedy}
\begin{algorithmic}
\REQUIRE A graph $G$ given with a list of edges and their weights 
\\ \hspace{1cm}$[W(\lbrace a,b \rbrace), \lbrace a,b \rbrace]$ which is sorted w.r.t weights. 
\ENSURE A Matching $M \subseteq E_G$.

\STATE $M \gets [\;]$
\STATE for $g$ in $G$ 
\STATE  \hspace{0.5cm}        if $g[2]$ has no intersection with (2nd component of) elements in $M$ then
\STATE  \hspace{1cm}             $\texttt{append}(M,g)$ 
\STATE return {$M$}
\end{algorithmic}
\end{algorithm}
\noindent
One can verify that the running time of the above algorithm is $O(m)$.

We complete this section by describing the construction of 
a bipartite graph corresponding to a given pair of polynomials. Here we use similar ideas which was used in~\cite{CorlessLeili2019Maple} except that this time we use the algorithm by~\cite{DuanPettie2014} to get a $\delta$-MWM.

Now we can form a bipartite weighted graph using the roots of the given polynomials $P$ and $Q$. Assume $R_P$ and $R_Q$ are sets of roots of $P$ and $Q$ respectively. Let $G^{\sigma}_{PQ} = (V_P,V_Q,E_{PQ}, W_{PQ})$ with

\begin{itemize}
\item $V_P= \texttt{ClusterRoots}(R_P, \sigma)$,
\item $V_Q= \texttt{ClusterRoots}(R_Q, \sigma)$,
\item $E^{\sigma}_{PQ} = \Big \lbrace \{r,s\}: \; r \in V_P,
s \in V_Q, \big\vert r-s \big\vert \leq \sigma \Big \rbrace$,
\item $W_{PQ}(\{r,s\}) = \min\{\text{multiplicity}_P(r),
\text{multiplicity}_Q(s)\}$.
\end{itemize}
A simple argument shows that we can form the corresponding graph to a
pair of polynomials by doing $n^2$ comparisons.\\

\section{Computing Approximate Polynomials and \GCD}\label{sec:apxpoly}
Assume $M$ is an MWM for $G^{\sigma}_{PQ}$. Each element in $M$ is an 
edge of the graph with a corresponding weight. We define
\begin{small}
\begin{equation}\label{eq:agcd}
G(x) = \prod_{\substack{\lbrace  a, b \rbrace \in M}} 
\left( x - \left( \frac{\mathrm{multiplicity}(a)\cdot a+\mathrm{multiplicity}(b) \cdot b}{\mathrm{multiplicity}(a)+\mathrm{multiplicity}(b)}\right)\right)^{W(\lbrace a,b \rbrace)}
\end{equation}
\end{small}
to be the approximate \GCD\  corresponding to $M$.


Using $\dm{t}$ (Equation~\eqref{pseudometric}) on the space of polynomials of degree $t$ we can find the 
approximate polynomials with desired properties. Suppose $G$ is the 
approximate \GCD\  as in Equation \eqref{eq:agcd}. Suppose $r_1, \ldots r_u$ are clustered roots of $P$ with
multiplicity $d_{i}$ which do not appear in the corresponding matching to 
$G$ and $r_{u+1}, \ldots, r_{\ell}$ with multiplicities $d_{u+1}, \ldots, d_{u+\ell}$ are roots of $P$ which appear in the corresponding matching to $G$. Similarly let $s_1, \ldots , s_v, s_{v+1}, \ldots s_{v+\ell}$ be the roots of $Q$ with multiplicity $d'_i$. Moreover suppose, $w_i = W(r_{u+i}, s_{v+i})$ for $1\leq i \leq \ell$.

We define
\begin{equation}
\tilde{P} = G(x)\cdot\prod_{\substack{1 \leq i \leq u}} \left( x- r_i\right)^{d_i} \cdot\prod_{\substack{u+1 \leq i \leq \ell}} \left( x- r_i\right)^{d_i-w_i}
\end{equation}

\begin{equation}
\tilde{Q} = G(x)\cdot\prod_{\substack{1 \leq i \leq v}} \left( x- s_i \right)^{d'_i} \cdot\prod_{\substack{v+1 \leq i \leq \ell}} \left( x- s_i \right)^{d'_i-w_i}
\end{equation}
Now it is not hard to see that there exists permutation of roots such that
$\dm{n}(P,\tilde{P}) \leq \sigma$ and $\dm{m}(Q, \tilde{Q})\leq \sigma$.

Since in each step of our algorithm we have access to roots
(and their multiplicities) of polynomials, presenting them in a Lagrange 
basis is straightforward.

%

\section{Numerical Results}
In this section we present a concrete example of computing an approximate \GCD\  by applying our algorithm to a pair of polynomials. In order to do so, we use a Maple implementation of our algorithm.  

Assume $P$ and $Q$ are given with the following 
data (in a Lagrange basis):
\begin{equation*}
\begin{split}
P_x = &[ 4.586334585, 5.161255391, 2.323567403, 1.809094426,\\ 
	& \; 1.471852626, 4.427838553, 2.275731771, 1.020544909]\\
P_y = & [ 787.1243900, 3285.933680, 0.01345240, 0.00001680,\\ 		& \; 0.00880350, 499.6500860, 0.01132600, 0.12249270]\\
Q_x = & [ 2.812852786, 1.746745227, 2.296707006, 2.359573808,\\ 
	& \; 4.747053250, 1.640439652, 5.832623175]\\
Q_y = &[- 0.0095256, 0.0171306, 0.2253058, 0.2359018,\\ 			 	& \; 426.4319036, 0.0041690, 3314.165173]
\end{split}
\end{equation*}
where $P_y[i]$'s are values of $P$ on $P_x[i]$'s and $Q_y[i]$'s are values of $Q$ on $Q_x[i]$'s. Since $P_x$ is of length 8 and 
$Q_x$ is of length 7, polynomials $P$ and $Q$ are respectively of degrees 7 and 6. 

 Asking Maple for generalized eigenvalue of the pair $(C_{P,0}, C_{P,1})$ gives
 $$  \left[ \begin {array}{c} 
  2.80+ 0.0\,i\\  2.60+ 0.0\,i
\\  1.90+ 0.0\,i\\  1.85+ 0.0\,i
\\  1.75+ 0.0\,i\\  1.70+ 0.0\,i
\\  1.0+ 0.0\,i
\end {array} \right] 
 $$
as roots of $P$ and the corresponding pair to $Q(x)$ gives
$$
\left[ \begin {array}{c} 
 3.0+ 0.0\,i\\  2.80+ 0.0\,i
\\  1.31+ 0.0\,i\\  1.33+ 0.0\,i
\\  1.45+ 0.0\,i\\  1.50+ 0.0\,i
\end {array} \right] 
$$
The second step is to cluster the set of roots. Before clustering we add multiplicities of the roots to the list of roots.
$$R_P = [[1, 1], [1.7, 1], [1.75, 1], [1.85, 1], [1.9, 1], [2.6, 1], [2.8, 1]],$$
$$R_Q = [[1.31, 1], [1.33, 1], [1.45, 1], [1.50, 1], [2.8, 1], [3, 1]].$$
By applying \texttt{ClusterRoots} with $\sigma = 0.5$ we get $$R_P = [[1, 1], [1.825000000, 4], [2.700000000, 2]],$$
$$R_Q=[[1.442500000, 4], [2.900000000, 2]]$$
The clustered roots of $P$ and $Q$ will give us the following bipartite graph

\begin{figure}[H]
	\includegraphics[width=5cm]{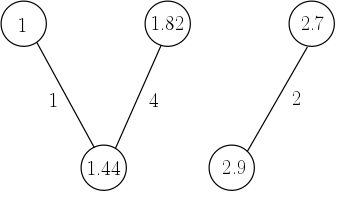}
	\centering
	\caption{The bipartite graph corresponding to $P$ and $Q$}.
\end{figure}
Note that there are two MCM's given by 
\begin{figure}[H]
    \centering
    \includegraphics[width=5cm]{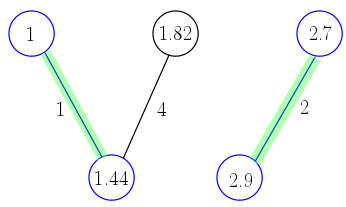}
    \caption{A Maximum Cardinality Matching for $G$, which is not a Maximum Weight Matching}\label{fig:MWM1}
\end{figure}

\begin{figure}[H]
    \centering
    \includegraphics[width=5cm]{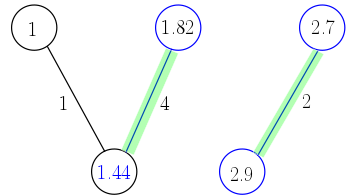}
    \caption{A Maximum Cardinality Matching for $G$, which is a Maximum Weight Matching as well.}\label{fig:MWM2}
\end{figure}
Although Figure~\ref{fig:MWM1} is a MCM, it is not a maximum weight matching. This
leads to a smaller degree in our approximate \GCD. However, Figure~\ref{fig:MWM2} 
is a maximum weight matching which gives us a higher degree \GCD.

Finally we need to find $G(x)$, $\tilde{P}$ and $\tilde{Q}$.\\
Since the MWM given in Figure~\ref{fig:MWM2} contains two edges we have $\ell = 2$ and
\begin{align*}
G(x) = & (x - \dfrac{1.82+1.44}{2})^4\cdot(x- \dfrac{2.7+2.9}{2})^2 \\
= & (x-1.63)^4\cdot(x-2.8)^2
\end{align*}
which will be presented by $G_x=[1.63,1.63,1.63,1.63,2.8,2.8, 0]$ and $G_y = [0,0,0,0,0,0,38.98397700]$.
In order to get $\tilde{P}$ and $\tilde{Q}$, we write:

\begin{align*}
&r_1 = 1, r_{1+1} = 1.82, r_{1+2} = 2.7\\
&d_1 = 1, d_{1+1} = 4, d_{1+2} = 2\\
&s_{0+1} = 1.44, s_{0+2} = 2.9\\
&d'_{0+1} = 4, d'_{0+2} = 2\\
&w_1 = 4, w_2 = 2\\
\end{align*}
This gives us 
\begin{align*}
&\tilde{P}(x) = G(x)\cdot (x-1)\\
&\tilde{Q}(x) = G(x)\\
\end{align*}
Also we can write 
\begin{align*}
&\tilde{P}_x= [1.63,1.63,1.63,1.63,2.8,2.8,1,0]\\
&\tilde{P}_y= [0,0,0,0,0,0,0,38.98397700]\\
\end{align*}

\section{Concluding Remarks\label{sec:concl}}
This paper solves the approximate \GCD\ problem in Lagrange bases. Our 
algorithm follows the ideas in~\cite{CorlessLeili2019Maple} and~\cite{Pan2001}. The results of this work are relevant to approximate \GCD\  for a specific pseudometric, namely the root pseudometric.

A possible extension of this work may be approximate \GCD\  in Lagrange bases
with respect to a different metric (or pseudometric). Moreover, the main
algorithm of this paper can
be trivially extended to Hermite basis. The only difference will be the rootfinding step. In order to compute the roots one can use the companion 
pencil presented in~\cite{CorlessChanLeili2018}.

One may wonder why clustering is used in $p(z)$ and separately in $q(z)$ to identify common roots, instead of simply considering the roots of $p(z)$ and $q(z)$ together and then clustering them (and thus identifying common zeros and hence the \texttt{GCD}).  We feel that a significant feature of approximate multiple zeros arising from a \textsl{single} polynomial is \textsl{symmetry}: computing $p(z)(z-r)^m + \varepsilon$ results in a near-circular symmetric cluster near $z=r$.  There is no reason to believe that the perturbation in $q(z)(z-r)^\ell+\delta$ would even be the same magnitude, and in any case symmetry would be lost.  We admit that this is a heuristic, and does not account for \textsl{structured} perturbations $p(z)(z-r)^n + \varepsilon r(z)$; nonetheless we think this is a valuable heuristic and it has worked well for us in the examples that we have tried.

An interesting example, which we will describe in a future paper, is to cluster the roots of a Mandelbrot polynomial, defined as $p_0(z)=0$ and $p_{n+1}(z) = zp_n^2(z) + 1$.  These are studied by homotopy methods in~\cite{Chan2016} using the homotopy $f(z,t) = zp_n^2(z) + t$ which starts with double roots at $t=0$: if one looks at the roots of $p_{n+1}(z)$ and clusters into double roots, do we get approximations of roots of $p_n(z)$?  The answer is no, or not often, and then only in regions of sensitivity.  But this example is, as mentioned in the previous paragraph, structured.  We need to perform more experiments on this example before commenting further.

Finally, as it is described in this current paper, the clustering heuristic by distance, multiplicity, and symmetry does not take into any account that the points are actually supposed to be roots of a polynomial.  Because it is based on the approximation
\begin{equation}
    p(r+\Delta r) \approx \frac{p^{(m)}(r)}{m!}\Delta r^m
\end{equation}
for an $m$-fold zero, we should in practice post-process the clustering results, estimate the $m$th derivative at a clustered root, and re-do the clustering with the ``fuzz factors'' of the heuristic replaced by the estimates resulting from the post-processing.  This would increase our confidence that we had found an actual multiple root.  Indeed one could imagine an iterative algorithm being based on this refinement step.  We leave this to future work.


\ifCLASSOPTIONcompsoc
  \section*{Acknowledgments}
\else
  \section*{Acknowledgment}
\fi
This work was supported by NSERC and the Ontario Research Centre for Computer Algebra.  RMC thanks Laureano Gonzalez-Vega for a conversation in 1997 about the potential value of symmetry in clustering.



%
\bibliographystyle{plain}
\bibliography{bib}

\end{document}